\documentclass[reqno,12pt]{amsart}

\usepackage{amsmath}
\usepackage{amssymb}
\usepackage{amsthm}
\usepackage[english]{babel}
\usepackage{graphics, xcolor}
\newtheorem{theorem}{Theorem}

\newtheorem{lemma}{Lemma}

\newtheorem{corollary}{Corollary}
\newcommand{\beqa}{\begin{eqnarray}}
\newcommand{\beqan}{\begin{eqnarray*}}
\newcommand{\eeqa}{\end{eqnarray}}
\newcommand{\eeqan}{\end{eqnarray*}}
\def\beq#1\eeq{\begin{equation}#1\end{equation}}

 \def\na{\,\, {\raise.4pt\hbox{$\shortmid$}}{\hskip-2.0pt\to}\, \, }

\def\={\overset{ \text{\rm def} }=}

\newcommand{\bl}{}

\newtheorem{remark}{Remark}

\begin{document}

\title[Concentration Functions]
{Improved applications of Arak's inequalities\\ to  the
Littlewood--Offord problem}

\author[F.~G\"otze]{Friedrich G\"otze}
\author[A.Yu. Zaitsev]{Andrei Yu. Zaitsev}

\email{goetze@math.uni-bielefeld.de}
\address{Fakult\"at f\"ur Mathematik,\newline\indent
Universit\"at Bielefeld, Postfach 100131,\newline\indent Bielefeld, D-33501,
Germany\bigskip}
\email{zaitsev@pdmi.ras.ru}
\address{St.~Petersburg Department of Steklov Mathematical Institute
\newline\indent
Fontanka 27,
\newline\indent
St.~Petersburg, 191023, Russia
\newline\indent
and
\newline\indent
St. Petersburg State University,
\newline\indent
Universitetskaya nab. 7/9,
\newline\indent
St. Petersburg, 199034, Russia}

\begin{abstract}
Let  $X_1,\ldots,X_n$ be independent identically distributed random variables.
In this paper we study the behavior of concentration functions of
 weighted sums $\sum_{k=1}^{n}X_ka_k $ with respect to the
arithmetic structure of coefficients~$a_k$ in the context of the
Littlewood--Offord problem. In our recent
 papers, we
discussed the relations between the inverse principles stated by
Nguyen, Tao and Vu and similar principles formulated  by Arak in
his papers from the 1980's. In this paper, we will
derive some improved (more general and more precise) consequences
 of Arak's inequalities
applying our new bound in the  Littlewood--Offord problem.
Moreover, we also obtain an improvement of the estimates used in Rudelson and Vershynin's least common denominator method.
\end{abstract}

\thanks{The authors were supported by SFB 1283/2 2021 –-- 317210226 and by the RFBR-DFG grant
 20-51-12004.}



\maketitle

\section{Introduction}

The aim of the present work is to provide a supplement to a paper of the authors~\cite{GZ18}. Using our recent results~\cite{GZ22}, we show that the dependence of constants on the involved distributions may be specified in more detail.

The
 concentration function of an $\mathbf{R}^d$-{valued}
vector $Y$ with distribution $F=\mathcal L(Y)$ is defined as
\begin{equation}
Q(F,\tau)=\sup_{x\in\mathbf{R}^d}\mathbf{P}(Y\in x+ \tau B), \quad
\tau\geq0, \nonumber
\end{equation}
where $B=\{x\in\mathbf{R}^d:\|x\|\leq 1/2\}$ 
denotes the centered
Euclidean ball of radius 1/2.

Let $X,X_1,\ldots,X_n$ be independent identically distributed
(i.i.d.{}) random variables. Let $a=(a_1,\ldots,a_n)\ne 0$, where
$a_k=(a_{k1},\ldots,a_{kd})\in \mathbf{R}^d$, $k=1,\ldots, n$.
Starting with seminal papers of Littlewood and Offord \cite{LO}
and Erd\"os~\cite{Erd}, the behavior of the concentration
functions of the weighted sums $S_a=\sum\limits_{k=1}^{n} X_k a_k$
has been intensively  studied. Denote by $F_a=\mathcal L(S_a)$ the
distribution of~$S_a$. We refer to~\cite{EGZ2} and \cite{Nguyen and Vu13} for a discussion of
the history of the problem (see also, for instance, Friedland and
Sodin \cite{Fried:Sod:2007}, Nguyen and Vu~\cite{Nguyen and Vu},
Rudelson and Vershynin \cite{Rud:Ver:2008}, \cite{Rud:Ver:2009},
Tao and Vu \cite{Tao and Vu}, \cite{Tao:Vu:2009:Bull},
 Vershynin \cite{Ver:2011}, Tikhomirov  \cite{T}, Livshyts,  Tikhomirov and Vershynin \cite{LTV}, Campos et al.
 \cite{CJMS}, \cite{CJMS22}).

The recent achievements in estimating the probabilities of singularity of random matrices  \cite{CJMS},  \cite{CJMS22}, \cite{LTV}, \cite{T},
 were based on the Rudelson and Vershynin \cite{Rud:Ver:2008}, \cite{Rud:Ver:2009}, \cite{Ver:2011} method of {\it least common denominator.}
Note that the results of \cite{Fried:Sod:2007}, \cite{Rud:Ver:2009}, \cite{Ver:2011} (concerning the
Littlewood--Offord problem) were improved and refined in \cite{Eliseeva}--\cite{EGZ}.

A few years ago, Tao and Vu \cite{Tao and Vu} and Nguyen and Vu
\cite{Nguyen and Vu} proposed the so-called {\it inverse principles} in
the Littlewood--Offord problem. In the
 papers  \cite{EGZ2} and~\cite{GZ18}, we
discussed the relations between these inverse principles and
similar principles formulated  by Arak
(see \cite{Arak0}--\cite{Arak and Zaitsev}) in his papers from the 1980's. In the one-dimensional
case, Arak has found a connection of the concentration function of
the sum with the arithmetic structure of supports of distributions
of independent random variables for {\it arbitrary}\/
distributions of summands. Using these results, he has solved in \cite{Arak} an old problem
stated by Kolmogorov~\cite{K}.

In \cite{GZ18}, we have shown that a consequence
of Arak's inequalities provides
results in the Littlewood--Offord problem of greater generality and improved precision compared to
those proved in \cite{EGZ2}. In the present paper, using our recent results~\cite{GZ22}, we show that the dependence of constants on the underlying distributions $\mathcal L(X)$  may be specified more precisely.  Moreover, in Section~\ref{s3}, we obtain an improvement of bounds obtained in the context of the method of the least common denominator.

Let us first introduce  the necessary notation.  Below
${\mathbf N}$ and ${\mathbf N}_0$ {will denote} the sets of all
positive and non-negative integers respectively. The symbol $c$
will be used for absolute positive constants.
 Note that $c$ can be different in different (or even in the same) formulas.
We will write $A\ll B$ if $A\leq c B$. {Furthermore,  we
{will} use the notation} $A\asymp B$ if $A\ll B$ and $B\ll A$.
If the corresponding constant depends on, say, $r$, we write
$A\ll_r B$ and $A\asymp_r B$.  If $\xi=(\xi_1,\ldots,\xi_d)$ is a
random vector with distribution~$F=\mathcal L(\xi)$, we denote
$F^{(j)}=\mathcal L(\xi_j)$, $j=1,\ldots,d$. Let
$\widehat F(t)={\mathbf E}\,\exp\big(\,i\,\langle
t,\xi\rangle\big)$, \,$t\in\mathbf R^d$, \, be the characteristic
function of the distribution~$F$. Here $\langle\, \cdot \,,\,
\cdot \,\rangle$ is the inner product in \,$\mathbf R^d$.

For~${x=(x_1,\dots,x_d )\in\mathbf R^d}$ we denote
$\|x\|^2=\langle x,x\rangle= x_1^2+\dots +x_d^2$ and $|x|=
\max_j|x_j|$. Let \,$E_y$ \,be the distribution concentrated at a
point ${y\in \mathbf{R}^d} $. We denote by $[K]_\tau$ the closed
$\tau$-neighborhood of a set $K$ in the norm
${|\,\cdot\,|}$. Products and powers of measures will be understood
in the sense of convolution. Thus, we write $W^k$ for the $k$-fold
convolution of a measure~$W$. For
a distribution~$W$, we will also denote \begin{equation}  e(\alpha\, W)=e^{-\alpha}\sum_{k=0}^\infty
 \frac{\alpha^k\,W^k}{k!}, \quad \alpha>0.  \label{00}\end{equation}
It is well-known that the distribution $D=e(\alpha\, W)$ is infinitely divisible with the L\'evy
 spectral measure $\alpha\, W$ and $\widehat D(t)=\exp\big(\alpha\,\big(\widehat W(t)-1\big)\big)$.
If a distribution~$D$ is
infinitely divisible, $D^\lambda$, $\lambda\ge0$, is the
infinitely divisible distribution with characteristic function
$\widehat D^\lambda(t)$. Note that $(e(\alpha\, W))^\lambda=e(\alpha\,\lambda\, W)$.

For a finite set $K$, we denote by $|K|$
the number of elements~$x\in K$. \bl {The symbol $\times$} is used
for the direct product of sets. We write $O(\,\cdot\,)$ if the
involved constants depend on the parameters named ``constants'' in
the formulations, but not on~$n$.

The elementary properties of concentration functions are well
studied (see, for instance,
 \cite{Arak and Zaitsev, Hengartner and Theodorescu,
Petrov}). In particular, it is clear that
\begin{equation}\label{8j}Q(F,\mu)\le \big(1+\lfloor
\mu/\lambda\rfloor\big)^d\,Q(F,\lambda), \quad\hbox{for any }
\mu,\lambda>0,\end{equation} where $\lfloor x\rfloor$ is the
largest integer~$k$ such that $k\le x$. Hence,
\begin{equation}\label{8a} Q(F,c\,\lambda)\asymp_d\,Q(F,\lambda).
\end{equation}
\medskip

Let  $r\in{\mathbf N}_0$, $m\in{\mathbf N}$ be fixed, let $h\in {\mathbf R}^r$
be an arbitrary $r$-dimensional vector,  and let $V$ be an
arbitrary closed symmetric convex subset of~${\mathbf R}^r$
containing not more than $m$ points with integer coordinates. We
define $\mathcal{K}_{r,m}$ as the collection of all sets of the
form
\begin{equation}
K=\big\{\langle{\nu}, h\rangle:{\nu}\in {\mathbf Z}^r\cap
V\big\}\subset{\mathbf R} .\label{1s1}
\end{equation}
We call such sets CGAPs (convex generalized arithmetic
progressions, see \cite{Green}), by analogy with generalized arithmetic
progressions (GAPs)
involved in recent investigations of the Littlewood--Offord problem. The definition of GAPs is given below.
In the  case $r=0$ the class $\mathcal{K}_{r,m}=\mathcal{K}_{0,m}$
consists of the single set $\{0\}$ having zero as the unique element.

For any Borel measure $W$ on ${\mathbf R}$ and $\tau\ge0$ we define $\beta_{r,m}(W, \tau)$ by
the relation
\begin{equation}
\beta_{r,m}(W, \tau)=\inf_{K\in\mathcal{K}_{r,m} }W\{{\mathbf
R}\setminus[K]_\tau\}, \label{1s3}
\end{equation}
where $[K]_\tau$ is the closed
$\tau$-neighborhood of the set $K$.

The following definition is given in Tao and Vu \cite{TV08} (see
also  \cite{TV06}).

Let $r\in{\mathbf N}_0$ be a non-negative integer,
$L=(L_1,\ldots, L_r)$
 be a $r$-tuple of positive reals,
and  $g=(g_1,\ldots, g_r)$ be a $r$-tuple of elements of $\mathbf
R^d$. The triplet $ P = (L, g, r)$  is called symmetric
generalized arithmetic progression (GAP) in $\mathbf R^d$. Here
$r$ is the rank, $L_1,\ldots, L_r$ are the dimensions and
$g_1,\ldots, g_r$ are the generators of the GAP~$P$. We define the
image $\hbox {Image\,}(P) \subset\mathbf R^d$ of $P$ to be the set
$$\hbox {Image}(P) = \big\{ m_1g_1 + \cdots + m_rg_r:-L_j \le m_j \le L_j,
 \ m_j\in{\mathbf Z}, \ j=1, \ldots,  r\big\}.
$$
  For $t>0$ we denote the $t$-dilate $P^t$ of $P$ as the symmetric
GAP $ P^t = (tL, g, r)$ with
$$\hbox {Image\,}(P^t) = \big\{ m_1g_1 + \cdots + m_rg_r:-tL_j \le m_j \le tL_j,
 \ m_j\in{\mathbf Z},\ j=1, \ldots,  r\big\}.
$$
We define the size of $P$ to be $\hbox{size}(P) = \left|\,\hbox
{Image}(P)\right|$.

In fact, $\hbox {Image}(P)$ is the image of an integer box $$B =
\big\{(m_1,\ldots, m_r) \in{\mathbf Z}^r: -L_j \le m_j \le
L_j\big\}$$ under the linear map
$$
\Phi:(m_1,\ldots m_r) \in{\mathbf Z}^r \rightarrow  m_1g_1 +
\cdots + m_rg_r.
$$
We say that $P$ is {\it proper} if this map is one to one, or,
equivalently, if
\begin{equation}\label{6u2}\hbox{size}(P)=\prod_{j=1}^{r}\big(2\,\lfloor L_j\rfloor+1\big).
\end{equation}
For non-proper GAPs, we  have, of course,
\begin{equation}\label{6u277}\hbox{size}(P)<\prod_{j=1}^{r}\big(2\,\lfloor
L_j\rfloor+1\big).
\end{equation}

 In the  case $r=0$ the vectors $L$ and
$g$ have no elements and the image of the GAP~$P$ consists of the
unique zero vector $0\in\mathbf R^d$.

Recall that a convex body in the $r$-dimensional Euclidean space
$\mathbf R^r$ is a compact convex set with non-empty interior.
The following Lemma \ref{llmj} is contained in Theorem~1.6 of
\cite{TV08}.

\begin{lemma}\label{llmj}
 Let $V$ be a convex symmetric body in
$\mathbf R^r$, and let $\Lambda$ be a lattice in $\mathbf R^r$.
Then there exists a symmetric  proper GAP $P$ in
$\Lambda$ with rank $l\le r$ such that we have
\begin{equation}
\label{6s2}\hbox {\rm Image}(P) \subset V\cap\Lambda\subset\hbox
{\rm Image} \bigl(P^{(c_1r)^{3r/2}}\bigr)
\end{equation}with an
absolute constant~$c_1\ge1$.
\end{lemma}

\begin{corollary}\label{clmj}Under  the conditions of Lemma\/
$\ref{llmj}$, the inequality
 \begin{equation} \hbox{\rm size}(P^{(c_1r)^{3r/2}})\le\big(2\,(c_1\,r)^{3r/2}+1\big)^r\left|\,V\cap\Lambda\,\right|\label{7s2}
\end{equation} holds.
\end{corollary}

Estimating the concentration functions $Q(F_a, \tau)$ in the Littlewood--Offord
problem, it is possible to reduce the problem to the estimation of
concentration functions of some symmetric infinitely divisible
distributions. The corresponding statements are contained in Lemmas
\ref{lm42} and \ref{l429} below.

 Introduce the distribution
$H$ with the characteristic function
\begin{equation} \label{11}\widehat{H}(t)
=\exp\Big(-\frac{\,1\,}2\;\sum_{k=1}^{n}\big(1-\cos\left\langle \,
t,a_k\right\rangle \big)\Big), \quad t\in\mathbf R^d.
\end{equation}
It is clear that $H$ is a symmetric infinitely divisible
distribution.

Let $\widetilde{X}=X_1-X_2$ \,be the symmetrized random variable
where $X_1$ and $X_2$ are i.i.d.\ variables involved in the
definition of $S_a$. In the sequel we use the notation
$G=\mathcal{L}(\widetilde{X})$. For $\delta\ge0$, we denote
\begin{equation}\label{pp}p(\delta)= G\big\{\{z:|z| >
\delta\}\big\}.\end{equation}

\begin{lemma}\label{lm42}
For any $\varkappa,\tau>0$, we have
\begin{equation}\label{1155}
Q(F_a, \tau) \ll_d Q(H^{p(\tau/\varkappa)}, \varkappa).
\end{equation}
\end{lemma}

Lemma~\ref{lm42} was the main tool in our recent applications of Arak's inequalities to  the
Littlewood--Offord problem (see \cite{EGZ2} and~\cite{GZ18}).

We should note that
$H^{b}$, $b\ge0$, is a symmetric infinitely divisible
distribution with the L\'evy
 spectral measure $\frac{\,nb\,}2\;M^*$, where
 \begin{equation}\label{11556}
  M^*=\frac{\,1\,}{2n}\sum_{k=1}^{n}\big(E_{a_k}+E_{-a_k}\big).
 \end{equation}
 In fact, \begin{equation}\label{11557}H^{b}=e\Big(\frac{\,nb\,}2\;M^*\Big).
\end{equation}

 Arak \cite{Arak0,Arak} has found  connections of values of  the
concentration functions of infinitely divisible distributions  with the arithmetic structure of
supports of their spectral measures. Concerning Lemma~\ref{lm42}, these connections imply that if   the
concentration function of an infinitely divisible distribution is large, then
its spectral measure is concentrated near a CGAP $K$ having some simple arithmetic structure. According to Lemma~\ref{llmj} and Corollary~\ref{clmj}, this CGAP $K$ is contained in the image of a GAP having  rank and size which are comparable with rank and size of~$K$.
The proof of Lemma~\ref{lm42} is given
in~\cite{EGZ2}. It is rather elementary and is based on known
properties of concentration functions.

According to \eqref{8j}, Lemma~\ref{lm42} implies the following
Corollary~\ref{lm429}.

\begin{corollary}\label{lm429}For any $\varkappa,\tau, \delta>0$, we have
\begin{equation}\label{1166}
Q(F_a, \tau) \ll_d\big(1+\lfloor\varkappa/\delta\rfloor \big)^d\,
Q(H^{p(\tau/\varkappa)}, \delta).
\end{equation}
\end{corollary}

It is reasonable to be free in
the choice of $\delta$ in~\eqref{1166}. In our recent paper~\cite{GZ22}, a more precise
statement than Lemma \ref{lm42} is obtained, see Lemma \ref{l429} below. It provides useful
bounds if $p(\tau/\varkappa)$ is small, even if
$p(\tau/\varkappa)=0$.

 \begin{lemma}\label{l429}
For any\/ $\varkappa,\tau>0$, we have
\begin{equation}\label{c234}
Q(F_a,\tau) \ll_d
\lambda^{-1}\,Q(H^{\lambda},\varkappa),
\end{equation}where
\begin{equation}\label{cc239}
\lambda=\lambda_d(\tau/\varkappa)=\int\limits_{z\in\mathbf{R}}\big(1+\lfloor\tau(\varkappa
|z|)^{-1}\rfloor \big)^{-d}\,G\{dz\}.
\end{equation}
\end{lemma}

It is clear that $\lfloor\tau(\varkappa
|z|)^{-1}\rfloor=0$ if $|z|>\tau/\varkappa$. Therefore. $\lambda=\lambda_d(\tau/\varkappa)\ge p(\tau/\varkappa)$ and hence $Q(H^{\lambda},\varkappa)\le Q(H^{p(\tau/\varkappa)},\varkappa)$. Thus,
if
$\lambda\gg_d 1$ is essentially larger than $p(\tau/\varkappa)$, then inequality~\eqref{c234} of Lemma~\ref{l429} is stronger than inequality~\eqref{1155} of Lemma~\ref{lm42}. Therefore, using inequality~\eqref{c234} implies more precise bounds than those obtained earlier on the base of
inequality~\eqref{1155}.

Lemmas~\ref{lm42} and~\ref{l429} are quite general since $\varkappa,\tau>0$ are arbitrary. Moreover, in the right-hand sides of
inequalities~\eqref{1155} and~\eqref{c234} the dependence of $F_a$ on $\mathcal{L}(X)$ is reduced to the dependence on $p(\tau/\varkappa)$
or $\lambda_d(\tau/\varkappa)$ while the dependence  on $a$ is reduced to the dependence on $H$. Note that the distribution $H$ also could be written in the form $H=F_a=\mathcal{L}(S_a)$ with
  \begin{equation}\label{11557u}\mathcal{L}(X)=e\Big(\frac{\,1\,}4\;E_{1}\Big)*e\Big(\frac{\,1\,}4\;E_{-1}\Big)
\end{equation}
being the symmetrized Poisson distribution with parameter 1/4.

 \bigskip

Now we formulate Theorem~\ref{thm2} which is a
one-dimensional Arak type result, see~\cite{Arak, Arak and Zaitsev}.
It is a particular case of Theorem~4.3 of Chapter~II  in \cite{Arak and Zaitsev}.

\begin{theorem}\label{thm2} Let $D=e(\alpha\, W)$ be a one-dimensional
infinitely divisible distribution with characteristic function of
the form $\exp\big(\alpha\,(\widehat W(t)-1)\big)$, $t\in{\mathbf
R}$, where $\alpha>0$ and\/ $W$ is a probability distribution.
Let\/ $\tau>0$, \;$r\in{\mathbf N}_0$, $m\in{\mathbf N}$. Then
\begin{equation}
Q(D, \tau)\le
c_2^{r+1}\biggl(\frac{1}{m\sqrt{\alpha\,\beta_{r,m}(W, \tau)}}
+\frac{(r+1)^{5r/2}}{(\alpha\,\beta_{r,m}(W,
\tau))^{(r+1)/2}}\biggr), \label{1s4}
\end{equation}
where $c_2$ is an absolute constant.
\end{theorem}

Corollary~\ref{lm429} and
Theorem~\ref{thm2} imply the following Theorem~\ref{pthm7}.

 \begin{theorem}\label{pthm7} Let $\varkappa, \delta,\tau>0$,   $d=1$, \;$r\in{\mathbf N}_0$, $m\in{\mathbf N}$. Then
\begin{multline}
Q(F_a, \tau)\ll_d c_3^{r+1}\big(1+\lfloor\varkappa/\delta\rfloor
\big)\,\biggl(\frac{1}{m\sqrt{n\,p(\tau/\varkappa)\,\beta_{r,m}(M^*, \delta)}}\\
+\frac{(r+1)^{5r/2}}{(n\,p(\tau/\varkappa)\,\beta_{r,m}(M^*, \delta))^{(r+1)/2}}\biggr),
 \label{1sy4}
\end{multline}
where
 $M^*$ is defined in~\eqref{11556} and $c_3$ is an absolute constant.
\end{theorem}

In order to prove Theorem~{\ref{pthm7}}, it suffices to apply
Corollary~\ref{lm429} and Theorem~\ref{thm2},
using that $H^{p(\tau/\varkappa)}$ is an infinitely
divisible distribution whose L\'evy
 spectral measure is $n\,p(\tau/\varkappa)\,M^*/2$ (see~\eqref{11557}).
 Introduce as well $M=\frac 1{2n}\sum_{k=1}^{n}E_{a_k}$.
 It is obvious that $M\le M^*$ and $\beta_{r,m}(M, \delta)\le\beta_{r,m}(M^*, \delta)$.
 It is clear that the value of the measure $\sum_{k=1}^{n}E_{a_k}$ on a set~$K$ is the number of $a_k$ belonging to $K$.
 Therefore, roughly speaking, $\beta_{r,m}(M, \delta)$ is $\frac 1{2n}$ multiplied by the number of $a_k $ not  belonging to the $\delta$-neighborhood of a CGAP $K\in\mathcal{K}_{r,m}$.

Similarly,  Lemma~\ref{l429} and
Theorem~\ref{thm2} imply the following Theorem~\ref{Th3}.

 \begin{theorem}\label{Th3} Let $\varkappa, \delta,\tau>0$,
$d=1$, \;$r\in{\mathbf N}_0$, $m\in{\mathbf N}$. Assume that ${\lambda_1(\tau/\varkappa)\gg1}$. Then
\begin{multline}
Q(F_a, \tau)
\ll_d c_4^{r+1}\big(1+\lfloor\varkappa/\delta\rfloor
\big)\,\biggl(\frac{1}{m\sqrt{n\,\beta_{r,m}(M^*, \delta)}}
\\
+\frac{(r+1)^{5r/2}}{(n\,\beta_{r,m}(M^*, \delta))^{(r+1)/2}}\biggr),  \label{t1sy4}
\end{multline}
where
 $M^*$ is defined in~\eqref{11556} and $c_4$ is an absolute constant.
\end{theorem}

 Let\/  $d=1$,  and\/ $q=Q(F_a, \tau)$,   $\varkappa,\tau>0$. If $q$ is relatively  large (for example,  if $q\ge n^{-A}$ as in Tao and Vu \cite{Tao and Vu} and Nguyen and Vu
\cite{Nguyen and Vu}), then Theorem~{\ref{Th3}} implies that $n\,\beta_{r,m}(M^*, \delta)$ and hence $n\,\beta_{r,m}(M, \delta)$ are small.
 Therefore, recalling that $M=\frac 1{2n}\sum_{k=1}^{n}E_{a_k}$,  we see that the set of weights $\{a_k\}$ is well approximated by the $\delta$-neighborhood $[K]_\delta$ of a CGAP~$K\in \mathcal{K}_{r,m}$.

\medskip

Let  $r\in{\mathbf N}_0$, $s\in{\mathbf N}$ be fixed, let $h$
be an arbitrary $r$-dimensional vector,  and let $P$ be a symmetric
GAP  with image $\hbox {\rm Image}(P)\subset\mathbf Z^r$
containing not more than $s$ points with integer coordinates. We
define $\mathcal{P}_{r,s}$ as the collection of all sets of the
form
\begin{equation}
K=\big\{\langle{\nu}, h\rangle:{\nu}\in \hbox {\rm Image}(P)\big\}\subset{\mathbf R} .\label{1s12}
\end{equation}

For any Borel measure $W$ on ${\mathbf R}$ and $\tau\ge0$ we define $\gamma_{r,s}(W, \tau)$ by the relation
\begin{equation}
\gamma_{r,s}(W, \tau)=\inf_{K\in\mathcal{P}_{r,s} }W\{{\mathbf
R}\setminus[K]_\tau\}. \label{1s3}
\end{equation}

The main result of Zaitsev~\cite{Z16} is the following Theorem~\ref{thm52}. Theorem~\ref{thm52} is stated in terms of symmetric
GAPs. It follows from Theorem~\ref{thm2} and Lemma~\ref{llmj}.

\begin{theorem}\label{thm52} Let $D$ be a one-dimensional
infinitely divisible distribution with characteristic function of
the form $\exp\big(\alpha\,(\widehat W(t)-1)\big)$, $t\in{\mathbf
R}$, where $\alpha>0$ and\/ $W$ is a probability distribution.
Let\/ $\tau\ge0$, \;$r\in{\mathbf N}_0$, $s\in{\mathbf N}$. Then
\begin{equation}
Q(D, \tau)\le
c_5^{r+1}\biggl(\frac{(c_6\,r+1)^{3r^2/2}}{s\sqrt{\alpha\,\gamma_{r,s}(W, \tau)}}
+\frac{(r+1)^{5r/2}}{(\alpha\,\gamma_{r,s}(W,
\tau))^{(r+1)/2}}\biggr), \label{16s4}
\end{equation}
where $c_5, c_6$ are absolute constants.
\end{theorem}

Similarly,  Lemma~\ref{l429} and
Theorem~\ref{thm52} imply the following Theorem~\ref{tthm75}.

 \begin{theorem}\label{tthm75} Let $\varkappa, \delta,\tau>0$,
$d=1$, \;$r\in{\mathbf N}_0$, $s\in{\mathbf N}$. Assume that ${\lambda_1(\tau/\varkappa)\gg1}$. Then
\begin{multline}
Q(F_a, \tau)
\ll_d c_7^{r+1}\big(1+\lfloor\varkappa/\delta\rfloor
\big)\,\biggl(\frac{(c_8\,r+1)^{3r^2/2}}{s\sqrt{n\,\gamma_{r,s}(M^*, \delta)}}
\\
+\frac{(r+1)^{5r/2}}{(n\,\gamma_{r,s}(M^*, \delta))^{(r+1)/2}}\biggr),  \label{t1sy45}
\end{multline}
where
 $M^*$ is defined in~\eqref{11556} and $c_7, c_8$ are absolute constants.
\end{theorem}

\section{Inverse principles}\label{2}

In \cite{EGZ2} and
\cite{GZ18},
 we compared our results with
the results of Nguyen, Tao and Vu \cite{Nguyen and Vu, Nguyen and
Vu13, Tao and Vu, Tao and Vu3}.

 A few years ago, Tao and Vu
\cite{Tao and Vu} formulated in the discrete case (with
$\tau=0)$ the so-called {"inverse principle" in the  Littlewood--Offord problem}, stating that
\begin{multline}\hbox{\it A set $a=(a_1,\ldots,a_n)$ with large small
ball probability}\\ \hbox{\it  must have strong additive structure.}\end{multline} Here
``large small ball probability'' means that
$$Q(F_a,0)=\max_x{\mathbf P}\{S_a=x\}\ge n^{-A}$$ with some constant
$A>0$. ``Strong additive structure'' means that a large part of
vectors $a_1,\ldots,a_n$ belong to a GAP of bounded size.

 Nguyen and Vu \cite{Nguyen and Vu}
have extended {this inverse principle} to the continuous case (with strictly positive argument of the concentration function
$\tau_{n}>0)$ obtaining, in particular, the following result.

\begin{theorem}\label{tnv2}
Let $X$ be a real random variable satisfying condition
\begin{equation}
G\{\{x\in\mathbf{R}\colon C_1<|x| < C_2\}\}\ge C_3,
\label{eq16}
\end{equation}
with positive constants $C_1,C_2,C_3$. Let $0 < \varepsilon < 1$,
$A>0$ be constants and $\tau=\tau_{n} > 0$ be a parameter that may
depend on $n$. Suppose that $a=(a_1,\ldots,a_n) \in
({\mathbf{R}^d})^n$ is a multi-subset of ${\mathbf R}^d$ such that
  $q=Q(F_a, \tau)\ge n^{-A}$.
 Then,
for any number $n'$ between $n^\varepsilon$ and $n$,
 there exists a symmetric proper GAP $P$ with image
 $K$ such that

 $1$. At least $n-n'$ elements of $a$ are  $\tau$-close to $K$.

 $2$. $P$ has small rank $r=O(1)$, and small size
 \begin{equation}
|K|\le \max\big\{O\big(q^{-1}(n')^{-1/2}\big), 1\big\}.
\label{12sp}
\end{equation}
 \end{theorem}

Theorem~\ref{tnv2} says that a large part of
vectors $a_1,\ldots,a_n$ belongs to a $\tau$-neighborhood of the image of a GAP with bounded size.

Theorem~3 of \cite{GZ18} (which is Theorem~\ref{thm17a} below) has a similar formulation in the case $d=1$. It is derived using Theorem~{\ref{pthm7}} of the present paper and also include  positive parameters $\tau, \varkappa, \delta$ such that $\delta\le\min\{\varkappa, \tau\}$. The quantity $\delta$ is responsible for the size of the neighborhood of the image of a GAP  which covers  a large part of
vectors $a_1,\ldots,a_n$. The number $n'$ is again the number of $a_j$ which may be not approximated.
The condition $n^\varepsilon\le n'\le n$ is replaced, for an {\it  arbitrary} fixed~$r\in{\mathbf N}_0$, by
\eqref{2b1}.
The number $r$ is here the rank of the approximating GAP. Its size is estimated from above in
  \eqref{2s4478}.
Comparing \eqref{2s4478} with \eqref{12sp}, we see that in \eqref{2s4478} the dependence of constants on $\mathcal{L}(X)$ is given in an explicit form via the dependence on $p(\tau/\varkappa)$. The same may be said about condition~\eqref{2b1}.  This is connected with the use of Theorem~{\ref{pthm7}}. Our Theorem~\ref{Th3} is stronger than Theorem~{\ref{pthm7}} provided that ${\lambda_1(\tau/\varkappa)\gg1}$ since $p(\tau/\varkappa)\le \lambda_1(\tau/\varkappa)$.
Note also that, in the case where $\varkappa=\delta=\tau$ in Theorem~\ref{thm17a}, one can choose $r=O(1)$ so that conditions $q\ge n^{-A}$ and $n^\varepsilon\le n'\le n$  imply \eqref{2b1}, for $n$ large enough. Therefore, the one-dimensional version of  Theorem \ref{tnv2} follows from   Theorem~\ref{thm17a}.

Our main results are Theorems \ref{thm17} and
\ref{thm278}. They are improvements of Theorems 3 and 9 of~\cite{GZ18} which were proved via Lemma~\ref{lm42} of the present paper.
Now we can use more precise Lemma~\ref{l429} replacing, in the statements of \cite{GZ18},  $p(\tau/\varkappa)$  by $\lambda_1(\tau/\varkappa)$. Note that $p(\tau/\varkappa)$ can be small, even equal  to zero, while $\lambda_1(\tau/\varkappa)$ is always strictly positive if the distribution $G$ is non-degenerate.

 \begin{theorem}\label{thm17a}
Let\/  $d=1$,  \;$a=(a_1,\ldots,a_n) \in {\mathbf{R}^n}$, $\varkappa, \delta,\tau>0$, $\delta\le\min\{\varkappa, \tau\}$,
and\/ $q=Q(F_a, \tau)$. There exist positive
absolute constants $c_9$--$c_{12}$
 such that  for any fixed $r\in{\mathbf N}_0$,
 and any $n'\in \mathbf N$ satisfying the inequalities
\begin{equation}\big(\,2\,c_9^{r+1}\,(r+1)^{5r/2}\,\varkappa\big/q\,\delta\,\big)^{2/(r+1)}/p(\tau/\varkappa)
< n' \le n, \label{2b1}
\end{equation}
 there exist $m\in \mathbf N$ and CGAPs $K^*, K^{**}\subset\mathbf R$ having ranks
 $\le r$ and sizes~$\ll m$ and $\ll c(r)\,m$ respectively and such that

$1$. At least $n-2\,n'$ elements $a_{k}$ of\/ $a$ are
$\delta$-close to $K^*$, that is, $a_{k}\in[K^*]_{\delta}$ $($this
means that for these elements $a_{k}$ there exist $y_{k}\in K^*$
such that $\left|a_{k}-y_{k}\right|\le\delta)$.

$2$. The number $m$ satisfies the inequality
 \begin{equation}
 m\le\frac{2\,c_9^{r+1}\,\varkappa}{ q\,\delta\,\sqrt{p(\tau/\varkappa)\,
n'}}+1, \label{2s4478}
\end{equation}

$3$.  The set $K^*$ is contained in the image $\overline K$ of a
symmetric GAP $\overline P$ which has rank \,$\overline l\le r$,
 size~$\ll (c_{10}\,r)^{3r^2/2}m$ and generators $\overline g_j$,
$j=1,\ldots,\overline l$, satisfying inequality
 $\bigl|\overline g_j\bigr|\le 2\,r
\left\|a\right\|/\sqrt{n'}$.

$4$.  The set $K^*$ is contained in the image $\overline{\overline
K}$ of a proper symmetric GAP $\overline{\overline P}$ which has
rank \,$\overline{\overline l}\le r$ and size~$\ll
(c_{11}\,r)^{15r^2/2}m$.

$5$. At least $n-2\,n'$ elements of\/ $a$ are  $\delta$-close to
 $K^{**}$.

$6$.  The set $K^{**}$ is contained in the image $\widetilde K$ of
a proper symmetric GAP ${\widetilde P}$ which has rank
\,$\widetilde l\le r$, size~$\ll (c_{12}\,r)^{21r^2/2}m$ and
generators $\widetilde g_j$, $j=1,\ldots,\widetilde l$, satisfying
the inequality
 $\bigl|\widetilde g_j\bigr|\le 2\,r
\left\|a\right\|/\sqrt{n'}$.
\end{theorem}

 \begin{theorem}\label{thm17} The statement of Theorem~$\ref{thm17a}$ remains true  after replacing  $p(\tau/\varkappa)$  by $\lambda_1(\tau/\varkappa)$ provided that $\lambda_1(\tau/\varkappa)\gg1$.
\end{theorem}

The formulation of Theorem~\ref{thm17a} is rather cumbersome, but
this is the price for its generality. The statement  may be
simplified in particular cases, for example, for
$\varkappa=\delta$ or for $\varkappa=\tau$.

The assertion of Theorem~\ref{thm17a} is
non-trivial for each fixed $r$ starting with $r=0$. In this
case $m=1$ and Theorem~\ref{thm17a} gives a bound for the number
 of elements $a_{k}$ of~$a$ which are outside of the interval
$[-\delta,\,\delta]$.

Theorem~\ref{thm17a} and hence Theorem~\ref{thm17} are  formulated for one-dimensional $a_k$, $k=1,\ldots, n$.
  However, it may be shown that
Theorem~\ref{thm17a} provides sufficiently rich  arithmetic
properties for the set $a=(a_1,\ldots,a_n) \in {(\mathbf{R}^d)}^n$
in the multivariate case as well (see Theorem~\ref{thm27} below).
It suffices to apply Theorem~\ref{thm17a} to the distributions
$F_a^{(j)}$, $j=1,\ldots,d$,
of coordinates of the vector~$S_a$.\medskip

 Introduce the vectors $a^{(j)}=(a_{1j},\ldots,a_{nj})$,
$j=1,\ldots, d$. It is obvious that $F_a^{(j)}=F_{a^{(j)}}$.

\begin{theorem}\label{thm27}
Let\/  $d>1$,  $q_j=Q(F_a^{(j)}, \tau_j)$,
$\varkappa_j, \delta_j,\tau_j>0$, $\delta_j\le\min\{\varkappa_j, \tau_j\}$, $j=1,\ldots,d$. Below\/ $c_9$--$c_{12}$ are positive
absolute constants from Theorem~$\ref{thm17a}$. Suppose that\/
$a=(a_1,\ldots,a_n) \in ({\mathbf{R}^d})^n$ is a multi-subset of\/
${\mathbf R}^d$.
 Let\/   $r_j\in \mathbf N_0$,
 and\/  $n'_j\in \mathbf N$, $j=1,\ldots,d$, satisfy inequalities\/
 \begin{equation}\big(\,2\,c_{9}^{r_j+1}\,(r_j+1)^{5r_j/2}\,
 \varkappa_j\big/q_j\,\delta_j\,\big)^{2/(r_j+1)}\big/p(\tau_j/\varkappa_j)\le n_j' \le n,
\label{2s45}
\end{equation}
 Then, for each $j=1,\ldots,d$, there exist $m_j\in \mathbf
N$  and CGAPs $K_j^*, K_j^{**} \subset\mathbf R$ having ranks $\le
r_j$ and sizes~$\ll m_j$ and $\ll c(r_j)\,m_j$ respectively and
such that

$1$. At least $n-2\,n'_j$ elements $a_{kj}$ of\/ $a^{(j)}$ are
$\delta_j$-close to
 $K_j^*$, that is,
$a_{kj}\in[K_j^*]_{\delta_j}$
 $($this means that for these elements $a_{kj}$  there exist $y_{kj}\in K_j^*$ such
that $\left|a_{kj}-y_{kj}\right|\le\delta_j)$.

$2$. $m_j$ satisfies inequality $m_j\le w_j$, where
 \begin{equation}
 w_j=\frac{2\,c_{9}^{r_j+1}\,\varkappa_j}{ q_j\,\delta_j\,\sqrt{p(\tau_j/\varkappa_j)\,
n'_j}}+1,\label{2s498}
\end{equation}

$3$.  The set $K_j^*$ is contained in the image $\overline K_j$ of
a symmetric GAP $\overline P_j$ which has rank \,$\overline l_j\le
r_j$, size~$\ll (c_{10}\,r_j)^{3r_j^2/2}m_j$ and generators
$\overline g_p^{(j)}$, $p=1,\ldots,\overline l_j$, satisfying
inequality
 $\bigl|\overline g_p^{(j)}\bigr|\le 2\,r_j
 \left\|a^{(j)}\right\|/\sqrt{n'_j}$.

$4$.  The set $K_j^*$ is contained in the image
$\overline{\overline K}_j$ of a proper symmetric GAP
$\overline{\overline P}_j$ which has rank \,$\overline{\overline
l}_j\le r_j$ and size~$\ll (c_6\,r_j)^{15r_j^2/2}m_j$.

$5$.  At least $n-2\,n'_j$ elements of\/ $a^{(j)}$ are
$\delta_j$-close to
 $K_j^{**}$.

$6$.  The set $K_j^{**}$ is contained in the image $\widetilde
K_j$ of a proper symmetric GAP ${\widetilde P}_j$ which has rank
\,$\widetilde l_j\le r_j$, size~$\ll (c_{12}\,r_j)^{21r_j^2/2}m_j$
and generators $\widetilde g_p^{(j)}$, $p=1,\ldots,\widetilde
l_j$, satisfying inequality
 $\bigl|\widetilde g_p^{(j)}\bigr|\le 2\,r_j
\left\|a^{(j)}\right\|/\sqrt{n'_j}$.
\end{theorem}

The multidimensional version of  Theorem \ref{tnv2} follows from  Theorem \ref{thm27}. This is stated in \cite[Theorem~5]{GZ18}.
 We used that the condition $q\ge n^{-A}$
of Theorem~\ref{tnv2} implies the condition $Q(F_a^{(j)},
\tau)\ge n^{-A}$, $j=1,\ldots,d$, of \cite[Theorem~5]{GZ18}, since
$Q(F_a^{(j)}, \tau)\ge Q(F_a, \tau)$.\emph{}

In Theorem \ref{thm27}, the approximating GAPs for~${a}$ are direct products of GAPs which are approximating for~${a^{(j)}}$.

\begin{theorem}\label{thm27a}The statement of Theorem~$\ref{thm27}$ remains true  after replacing
$p(\tau_j/\varkappa_j)$  by $\lambda_1(\tau_j/\varkappa_j)$ provided that $\lambda_1(\tau_j/\varkappa_j)\gg1$, $j=1,\ldots,d$.
\end{theorem}

Theorems \ref{thm17a}--\ref{thm27a} have
non-asymptotic character.  They provide information about the
arithmetic structure of $a=(a_1,\ldots,a_n)$ without assumptions
like $q\ge n^{-A}$, which are
imposed in Theorem \ref{tnv2}. Theorems \ref{thm17a}--\ref{thm27a} are formulated for {\it fixed~$n$} and the dependence
of constants on parameters and on $\mathcal{L}(X)$ is given explicitly.

Theorem~\ref{thm19} below is \cite[Theorem~7]{GZ18}. It is  a  consequence of
Theorem~\ref{thm17a}.

 \begin{theorem}\label{thm19}
 Let\/  $\theta, A,\varepsilon_1,\varepsilon_2,\varepsilon_3,\varepsilon_4>0$ and\/ $B,D\ge0$
be constants, and\/ $\theta>D$. Let\/ $b_n,\varkappa_n, \delta_n,
\tau_n,\rho_n>0$, $n=1,2,\ldots$, be depending on $n$ non-random
parameters satisfying the relations \;$p(\tau_n/\varkappa_n)\ge
\varepsilon_3\,b_n^{-D}$, \;$
\varepsilon_4\,b_n^{-B}\leq\rho_n=\delta_n/\varkappa_n\leq1$,
$\delta_n\le \tau_n$, for all $n\in \mathbf
N$, and\/ $b_n\to\infty$ as $n\to\infty$.
 Suppose that $a=(a_1,\ldots,a_n) \in
{(\mathbf{R}^d)}^n$ is a multi-subset of\/~${\mathbf R}^d$ such
that $q_j=Q(F_a^{(j)}, \tau_n)\ge \varepsilon_1\,b_n^{-A}$,
\;$j=1,\ldots,d$, \;$n=1,2,\ldots$.
 Then, for any numbers $n'_j$ such that $\varepsilon_2\,b_n^\theta\leq n'_j\leq n$, $j=1,\ldots,d$,
 there exists a  proper symmetric GAP $P$ such that

$1$. At least $n-2\sum_{j=1}^dn'_j$ elements of\/ $a$ are
$\delta_n$-close to the image $K$ of the GAP $P$ in the norm
$|\,\cdot\,|$.

 $2$. $P$ has small rank $L=O(1)$, and small size
 \begin{equation}
|K|\le
\prod_{j=1}^d\max\Big\{O\Big(q_j^{-1}\,\rho_n^{-1}\,(n'_j\,p(\tau_n/\varkappa_n))^{-1/2}\Big),
1\Big\}. \label{n111sp}
\end{equation}
\end{theorem}
\medskip

Theorem~\ref{thm19} is more general than Theorem~\ref{tnv2},
where we restricted ourselves to the case $b_n=n$, $\varkappa_n=\tau_n=\delta_n$, {$n_j'=n'$} only.
Theorem~\ref{thm270} is an improvement of Theorem~$\ref{thm19}$.

\medskip
\begin{theorem}\label{thm270}The statement of Theorem~$\ref{thm19}$ remains true  after replacing
$p(\tau_n/\varkappa_n)$  by $\lambda_1(\tau_n/\varkappa_n)$ provided that $\lambda_1(\tau_n/\varkappa_n)\gg1$.
\end{theorem}

Similarly as in \cite{EGZ2, GZ18}, now we state analogs of Theorems~\ref{thm17a}--\ref{thm270}
for GAPs of logarithmic rank and with special dimensions, all equal to 1.
Theorems~\ref{tnv3}--\ref{thm5} below are Theorems~9--11 of \cite{GZ18}. They extend Theorems~5--7 of~\cite{EGZ2}.

\begin{theorem}\label{tnv3} Let\/  $d=1$,  \;$a=(a_1,\ldots,a_n) \in {\mathbf{R}^n}$, $\varkappa, \delta,\tau>0$, $\delta\le\min\{\varkappa, \tau\}$,
and\/ $q=Q(F_a, \tau)$.
Then there exists  a GAP $P$ of rank~$r\in\mathbf N$, of size~$\le3^r$,
with generators $g_k\in\mathbf{R}$, $k=1,\ldots,r$, and  such that its image $K\subset\mathbf{R}$ has the form
\begin{equation}
{K}=\Big\{\sum_{k=1}^r m_k\, g_k:m_k\in \{-1,0,1\}, \hbox{
for }k=1,\ldots,r\Big\}.\label{1s17}
\end{equation}
Moreover,  we have
\begin{equation}
r\ll\left|\log
q\right|+\log(\varkappa/\delta)+1, \label{n1ss65}
\end{equation}
and at least $n-n'$ elements of\/ $a$ are  $\delta$-close to
 $K$, where $n'\in\mathbf N$ and
\begin{equation}
n'\ll\big(p(\tau/\varkappa)\big)^{-1}
\bigl(\left|\log q\right|+\log(\varkappa/\delta)+1\bigr)^3.
\label{n1ss68d}
\end{equation}
\end{theorem}

\begin{theorem}\label{tnv4}
    Let\/  $d>1$,  $q_j=Q(F_a^{(j)}, \tau_j)$,
$\varkappa_j, \delta_j,\tau_j>0$, $\delta_j\le\min\{\varkappa_j, \tau_j\}$, $j=1,\ldots,d$.
     Then, for each $j=1,\ldots,d$,  there exists a GAP $P_j$
     of rank~${r_j\in\mathbf N}$, of size~$\le3^{r_j}$,
with generators $g_k^{(j)}\in\mathbf{R}$, $k=1,\ldots,r_j$, and
such that its image $K_j\subset\mathbf{R}$ has the form
\begin{equation}
{K_j}=\Big\{\sum_{k=1}^{r_j} m_k\, g_k^{(j)}:m_k\in \{-1,0,1\}, \hbox{
for }k=1,\ldots,r_j\Big\}.\label{1s174}
\end{equation}
Moreover,
\begin{equation}
r_j\ll\left|\log
q_j\right|+\log(\varkappa_j/\delta_j)+1, \label{n1ss654}
\end{equation}
and at least $n-n_j'$ elements of\/ $a^{(j)}$ are  $\delta_j$-close to
 $K_j$, where $n'_j\in\mathbf N$ satisfy
\begin{equation}
n'_j\ll\bigl(p(\tau_j/\varkappa_j)\bigr)^{-1}
\bigl(\left|\log q_j\right|+\log(\varkappa_j/\delta_j)+1\bigr)^3.
\label{n1ss68d4}
\end{equation}

Define $K=\times
_{j=1}^{d}K_j$. Then the set $K$ is the image of the $d$-dimensional GAP $P$ with rank
\begin{equation}
R=\sum_{j=1}^{d}r_j\ll\sum_{j=1}^{d}\bigl(\left|\log
q_j\right|+\log(\varkappa_j/\delta_j)+1\bigr), \label{n1ss658}
\end{equation}
and such that at least  $n-\sum_{j=1}^dn'_j$ elements of\/ $a$ belong to
the set $\times
_{j=1}^{d}[K_j]_{\delta_j}$. Here
\begin{equation}
\sum_{j=1}^dn'_j\ll\sum_{j=1}^{d}\bigl(p(\tau_j/\varkappa_j)\bigr)^{-1}
\bigl(\left|\log q_j\right|+\log(\varkappa_j/\delta_j)+1\bigr)^3.
\label{n1ss68d8}
\end{equation}

Furthermore, the set $K$  can be
represented as \begin{equation}
{K}=\Big\{\sum_{k=1}^{R} m_s\, g_s:m_s\in \{-1,0,1\}, \hbox{
for }s=1,\ldots,R\Big\}.\label{1s58}
\end{equation} Moreover, every vector $g_s\in
{\mathbf R}^d$, $s=1,\ldots,R$, has one non-zero coordinate only. Denote
$$
 s_0=0\quad\hbox{and}\quad s_j=\sum_{m=1}^{j}r_m, \quad j=1,\ldots,d.
$$
For $s_{j-1}<s\le s_j$, the vectors $g_s$ are non-zero in the
$j$-th coordinates only and these coordinates are equal to the elements of the
sequence $g_1^{(j)},\ldots, g_{r_j}^{(j)}$ from \eqref{1s174}.
\end{theorem}

\begin{theorem}\label{thm5}  Let\/  $A>0$ and\/ $B\ge0$
be constants. Let\/ $b_n,\varkappa_n, \delta_n,
\tau_n>0$, $n=1,2,\ldots$, be depending on $n$ non-random
parameters satisfying the relations  \;$
b_n^{-B}\leq\delta_n/\varkappa_n\leq1$,
$\delta_n\le \tau_n$, for all $n\in \mathbf
N$, and\/ $b_n\to\infty$ as $n\to\infty$.
 Let
 $q_j=Q(F_a^{(j)}, \tau_n)\ge b_n^{-A}$, for  $j=1,\ldots,d$.
 Then, for each $j=1,\ldots,d$,  there exists a GAP $P_j$
 of rank~${r_j\in\mathbf N}$, of size~$\le3^{r_j}$,
with generators $g_k^{(j)}\in\mathbf{R}$, $k=1,\ldots,r_j$, and
such that its image $K_j\subset\mathbf{R}$ has the form
\begin{equation}
{K_j}=\Big\{\sum_{k=1}^{r_j} m_k\, g_k^{(j)}:m_k\in \{-1,0,1\}, \hbox{
for }k=1,\ldots,r_j\Big\}.\label{1s19}
\end{equation} Moreover, the set $K=\times
_{j=1}^{d}K_j$ is the image of the $d$-dimensional GAP $P$ with rank
\begin{equation}
R=\sum_{j=1}^{d}r_j\ll d\,\big((A+B)\,\log b_n+1\big), \label{21st65}
\end{equation}
and such that at least  $n-n'$ elements of\/ $a$ belong to
the set $\times
_{j=1}^{d}[K_j]_{\delta_n}$. Here $n'\in\mathbf N$ and
\begin{equation}
n'\ll d \,\bigl(p(\tau_n/\varkappa_n)\bigr)^{-1}\big((A+B)\,\log
b_n+1\big)^3. \label{21st68d}
\end{equation}
Furthermore, the description of the set $K$ given at the end of the statement of
Theorem~$\ref{tnv4}$ remains true.\end{theorem}

In Theorems 5--7 of~\cite{EGZ2}, we obtained particular cases
of Theorems~\ref{tnv3}--\ref{thm5}, where
$b_n=n$ and $\tau=\varkappa$, $\tau_j=\varkappa_j$, $j=1,\dots,d$, or $\tau_n=\varkappa_n$, $n\in\mathbf N$.
Theorems~\ref{tnv3}\nobreak--\ref{thm5} were proved in~\cite{GZ18} using Lemma~\ref{lm42} and  Theorem~3.3 of Chapter~II  in \cite{Arak and Zaitsev} obtained by Arak~\cite{Arak0}. In Theorems~\ref{tnv3}--\ref{thm5} the approximating GAP may be non-proper.
 Replacing   Lemma~\ref{lm42}
by Lemma \ref{l429} in the proofs, we replace again $p(\,\cdot\,)$ by  $\lambda_1(\,\cdot\,)$ in the formulations  of Theorems~\ref{tnv3}--\ref{thm5} and obtain the following Theorems~\ref{thm278}--\ref{thm2700}.
\medskip
\begin{theorem}\label{thm278}The statement of Theorem~$\ref{tnv3}$ remains true  after replacing
$p(\tau/\varkappa)$  by $\lambda_1(\tau/\varkappa)$ provided that $\lambda_1(\tau/\varkappa)\gg1$.
\end{theorem}\medskip
\begin{theorem}\label{thm279}The statement of Theorem~$\ref{tnv4}$ remains true  after replacing
$p(\tau_j/\varkappa_j)$  by $\lambda_1(\tau_j/\varkappa_j)$ provided that $\lambda_1(\tau_j/\varkappa_j)\gg1$, $j=1,\dots,d$.
\end{theorem}\medskip
\begin{theorem}\label{thm2700}The statement of Theorem~$\ref{thm5}$ remains true  after replacing
$p(\tau_n/\varkappa_n)$  by $\lambda_1(\tau_n/\varkappa_n)$ provided that $\lambda_1(\tau_n/\varkappa_n)\gg1$.
\end{theorem}\medskip

\begin{remark}\rm
  In the statements of Theorems~\ref{thm278}--\ref{thm2700} we could replace  $p(\,\cdot\,)$ not by $\lambda_1(\,\cdot\,)$, but simply  by 1, provided that $\lambda_1(\,\cdot\,)\gg1$.
\end{remark}

\section{The method of the least common denominator}\label{s3}

In this section, we will continue the estimation of $Q(F_a, \tau)$, where $\tau>0$,  $F_a=\mathcal L(S_a)$, $S_a=\sum\limits_{k=1}^{n} X_k a_k$,
$X,X_1,\ldots,X_n$ are
i.i.d.\ random variables, and $a=(a_1,\ldots,a_n)\ne 0$, where
$a_k=(a_{k1},\ldots,a_{kd})\in \mathbf{R}^d$, $k=1,\ldots, n$.

 We will apply our Lemma~\ref{l429} to obtain more precise bounds in terms of the least common denominator method of Rudelson and Vershynin \cite{Rud:Ver:2008}, \cite{Rud:Ver:2009}, \cite{Ver:2011}.

Introduce the matrices
\begin{equation}\label{1155u}
\mathbb{A}=\sum_{k=1}^{n}\mathbb{A}_k,\quad\mathbb{A}_k=\left(\begin{matrix}a_{k1}^2&a_{k1}a_{k2}&\ldots&a_{k1}a_{kd}\\
a_{k2}a_{k1}&a_{k2}^2&\ldots&a_{k2}a_{kd}\\
\ldots&\ldots&\ldots&\ldots\\
a_{kd}a_{k1}&a_{kd}a_{k2}&\ldots&a_{kd}^2
\end{matrix}\right).
\end{equation}
We will use  the same letter $\mathbb{A}$ to denote the corresponding linear operator $\mathbb{A}:\mathbf{R}^d\to\mathbf{R}^d$,

We will use notation of Rudelson and Vershynin~\cite{Rud:Ver:2009}: for $t\in\mathbf{R}^d$
\begin{equation}  \label{667}
t\,\cdot\, a=(\left\langle \,
t,a_1\right\rangle, \ldots,\left\langle\,
t,a_n\right\rangle)\in\mathbf{R}^n.
\end{equation}
Thus,
\begin{equation}  \label{667oo}
\left\|t\cdot a\right\|^2=\sum_{k=1}^{n}\left\langle \,
t,a_k\right\rangle^2=\left\langle
\mathbb{A}t,t\right\rangle.
\end{equation}

\begin{theorem}\label{t19}Let $\alpha, b,D>0$, $0<\gamma<1$.
Let the distribution $H$ be defined in~\eqref{11}.
Assume that
\begin{equation}  \label{4bbn}\begin{split}
\Big(\sum_{k=1}^{n}\big(\left\langle \,
t,a_k\right\rangle-m_k \big)^2\Big)^{1/2}\ge\min\big\{\gamma\left\|t\cdot a\right\|,\alpha\big\}\hbox{ for all }m_1,\ldots, m_n\in \mathbf{Z}
, \
\\\hbox{and for all } t\in\mathbf{R}^d \hbox{ such that } \left\|t\right\|\le D.
\end{split}\end{equation}
Then
\begin{equation}\label{1ds}
Q(H^b, {1}/{D}) \ll_d \Big(\frac1{\gamma D\sqrt b}\Big)^d\frac1{\sqrt{\det{\mathbb{A}}}}+\exp(-4b\alpha^2).
\end{equation}
\end{theorem}

 It is clear that if
 \begin{multline} \label{4bb}
0<D\le D(a)=D_{\gamma, \alpha}(a)\\ =
\inf\Big\{\theta>0:\hbox{dist}(t\,\cdot\, a,\mathbf{Z}^n)<\min\big\{\gamma\left\|t\cdot a\right\|,\alpha\big\}
\\ \hbox{ for some }t\in\mathbf{R}^d\hbox{ such that }\|t\|=\theta \Big\},
\end{multline}
where $$\hbox{dist}(t\,\cdot\, a,\mathbf{Z}^n)= \min_{m \in \mathbf{Z}^n}\|\,t\,\cdot\, a - m\|
=\Big(\sum_{k=1}^{n} \min_{m_k \in \mathbf{Z}} \big(\left\langle \,
t,a_k\right\rangle-m_k \big)^2\Big)^{1/2} ,
 $$then  condition \eqref{4bbn} holds.

Rudelson and Vershynin~\cite{Rud:Ver:2009} called $ D(a)$ {\it the essential least
common denominator} of a vector $a$. There exist some similar but different definitions of least
common denominator (see \cite{Fried:Sod:2007}, \cite{Ver:2011}). The definition~\eqref{4bb} seems to be the most popular one.

Let us formulate a generalization of the classical Ess\'een
inequality \cite{Ess:1966} to the multivariate case
(\cite{Ess:1968}, see also \cite{Hengartner and Theodorescu}):

\begin{lemma}\label{lm3} Let $\tau>0$ and let
 $F$ be a $d$-dimensional probability distribution. Then
\begin{equation}
Q(F, \tau)\ll_d \tau^d\int_{\left\|t\right\|\le1/\tau}|\widehat{F}(t)| \,dt.
\label{4s4d}
\end{equation}
\end{lemma}

Lemma \ref{lm3} says that $Q(F, \tau)$ is estimated by the mean value of $|\widehat{F}(t)|$ on the ball $\{t\in\mathbf{R}^d:\left\|t\right\|\le1/\tau\}$.

\emph{Proof of Theorem\/ $\ref{t19}$.} Using Lemma \ref{lm3} and relation \eqref{8a}, we get
\begin{eqnarray*}
Q(H^{b}, 1/D)&\ll_d& D^{-d}\int_{\left\|t\right\|\le 2\pi D}|\widehat{H}(t)|^b \,dt\\
&=&D^{-d}\int_{\left\|t\right\|\le 2\pi D}\exp\Big(-\frac{\,b\,}2\;\sum_{k=1}^{n}\big(1-\cos\left\langle \,
t,a_k\right\rangle \big)\Big) \,dt
\\ .
\label{4s4d}
\end{eqnarray*}

It is easy to show
that  $1-\cos x \geq 2x^2/\pi^2$, for~${|x|\leq\pi}$. For
arbitrary~$x$, this implies that
$$1-\cos x\geq 2\,\pi^{-2}
\min_{m\in \mathbf{Z}}|\,x-2\pi m|^2.$$ Substituting this
inequality into \eqref{11}, we obtain
\begin{eqnarray}
\widehat{H}(t)&\leq&\exp \Big(-\cfrac{1}{\pi^{2}} \;\sum_{k=1}^{n}\min_{m_k \in
\mathbf{Z}}\big| \left\langle \,
t,a_k\right\rangle -2 \pi m_k\big|^2\Big)\nonumber \\
&=&\exp\Big(- \cfrac{1}{\pi^{2}}\;\big(\hbox{dist}(t\cdot a,2\pi\mathbf{Z}^n)\big)^2\Big)\nonumber \\
&=&\exp\Big(- 4\;\big(\hbox{dist}(t/2\pi\cdot a,\mathbf{Z}^n)\big)^2\Big).\label{7b}
\end{eqnarray}

Hence, for $\left\|t\right\|\le 2\pi D$, we have
\begin{eqnarray}
\widehat{H}(t)&\leq&
\exp\Big(- 4\;\big(\min\big\{\gamma\left\|t/2\pi\cdot a\right\|,\alpha\big\}\big)^2\Big).\label{7boo}
\end{eqnarray}
Therefore, we have
\begin{eqnarray*}
Q(H^{b}, 1/D)&\ll_d& D^{-d}\int_{\left\|t\right\|\le 2\pi D}\exp\Big(- 4\,b\gamma^2\;\left\|t/2\pi\cdot a\right\|^2\Big) \,dt\\
&&\quad\quad\quad\qquad+\ D^{-d}\int_{\left\|t\right\|\le 2\pi D}\exp\big(- 4\,b\alpha^2\big) \,dt\\
&\ll_d& D^{-d}\int_{\left\|t\right\|\le 2\pi D}\exp\Big(- b\gamma^2\pi^{-2}\;\left\langle
\mathbb{A}t,t\right\rangle\Big) \,dt\\
&&\quad\quad\quad\qquad+\ \exp\big(- 4\,b\alpha^2\big)\\
&\ll_d& \Big(\frac1{\gamma D\sqrt b}\Big)^d\frac1{\sqrt{\det{\mathbb{A}}}}+\exp(-4b\alpha^2),
\label{4s4doo}
\end{eqnarray*}
proving Theorem \ref{t19}. $\square$

Applying Lemma~\ref{l429} with $\varkappa=1/D$ and Theorem~\ref{t19} with $b=\lambda_d(\tau D)$,  where $\lambda_d(\,\cdot\,)$
is defined in \eqref{cc239}, we obtain the following
Theorem \ref{t20}.

\begin{theorem}\label{t20}Let the conditions
of Theorem\/ $\ref{t19}$ be satisfied. Then, for any $\tau>0$,
\begin{equation}\label{1dsoo}
Q(F_a, \tau) \ll_d \frac1{\lambda_d(\tau D)}\bigg(\Big(\frac1{\gamma D\sqrt {\lambda_d(\tau D)}}\Big)^d\frac1{\sqrt{\det{\mathbb{A}}}}+\exp(-4\,\lambda_d(\tau D)\,\alpha^2)\bigg).
\end{equation}
\end{theorem}

Applying Lemma~\ref{lm42} with $\varkappa=1/D$ and Theorem~\ref{t19} with $b=p(\tau D)$,  where $p(\,\cdot\,)$
is defined in \eqref{pp}, we obtain the following
Theorem \ref{t21}.

\begin{theorem}\label{t21}Let the conditions
of Theorem\/ $\ref{t19}$ be satisfied. Then, for any $\tau>0$,
\begin{equation}\label{1dsoo}
Q(F_a, \tau) \ll_d \Big(\frac1{\gamma D\sqrt {p(\tau D)}}\Big)^d\frac1{\sqrt{\det{\mathbb{A}}}}+\exp(-4\,p(\tau D)\,\alpha^2).
\end{equation}
\end{theorem}

Theorem \ref{t21} is more precise and general than \cite[Theorem 3.3]{Rud:Ver:2009}.

Denote
\begin{equation} \label{0}M(\tau)=\tau^{-2}\int_{|x|\leq\tau}x^2
\,G\{dx\}+\int_{|x|>\tau}G\{dx\}=\mathbf{E}
\min\big\{{\widetilde{X}^2}/{\tau^2},1\big\}, \quad \tau>0.\end{equation}

The following Theorem \ref{t22} was obtained in \cite{Eliseeva}. Its one-dimensional version can be found in \cite{Eliseeva and Zaitsev}.

\begin{theorem}\label{t22}Let the conditions
of Theorem\/ $\ref{t19}$ be satisfied. Then, for any $\tau>0$,
\begin{equation}\label{1dsoo1}
Q(F_a, \tau) \ll_d \Big(\frac1{\gamma D\sqrt {M(\tau D)}}\Big)^d\frac1{\sqrt{\det{\mathbb{A}}}}+\exp(-c\,M(\tau D)\,\alpha^2).
\end{equation}
\end{theorem}

It is easy to see that $M(\tau D)\ge p(\tau D)$. Therefore, Theorem \ref{t22} is stronger than Theorem \ref{t21}.

\end{document}